\documentclass[a4paper, 11pt]{article}
\usepackage{bbm}
\usepackage{amsfonts}
\usepackage{amsmath}
\usepackage{amsthm,mathrsfs}
\usepackage{amsfonts,amssymb}
\usepackage{enumerate}
\usepackage{setspace}
\usepackage{authblk,titletoc}
\allowdisplaybreaks[4]

\newtheorem{theorem}{Theorem}[section]
\newtheorem{proposition}[theorem]{Proposition}

\newtheorem{lemma}[theorem]{Lemma}
\newtheorem{corollary}[theorem]{Corollary}

\newtheorem{remark}[theorem]{Remark}

\numberwithin{equation}{section}

\def\az{\alpha}

\def\A{\mathcal{A}}

\def\O{\mathcal{O}}

\def\ep{\epsilon}

\def\rr{{\mathbb R}}
\def\rn{{{\rr}^n}}

\def\mm{{\mathbb M}}

\def\az{{\alpha}}

\begin{document}
	\title{{Asymptotic expansion at infinity of solutions to Monge-Amp\`{e}re equation with $C^\alpha$ right term} }
	
	\author[1]{ Shuai Qi }
	\author[,1]{Jiguang Bao\thanks{Corresponding author.	}}
	\affil[1]{School of Mathematical Sciences, Beijing Normal University, 
		Laboratory of Mathematics and Complex Systems, Ministry of Education, 	
		Beijing, 100875, China}
	\date{}
	\maketitle
	
	\renewcommand\thefootnote{}
	\footnote{{\it E-mail address:} qshuai@pku.edu.cn(S. Qi), jgbao@bnu.edu.cn(J. Bao). }
	\setcounter{footnote}{0}
	\renewcommand{\thefootnote}{\fnsymbol{footnote}}
	\vspace{-1cm}
	
	{\bf Abstract:} We develop a non-local method to establish the asymptotic expansion at infinity of solutions to 
	Monge-Amp\`{e}re equation $\det(D^2v)=f$ on $\rn$, where $f$ is a perturbation of $1$ and is only assumed to be H\"{o}lder 
	continuous outside a bounded subset of $\rn$, compared to the previous work that $f$ is at least $C^2$.


	\abovedisplayskip=7pt
	\abovedisplayshortskip=7pt
	\belowdisplayskip=7pt
	\belowdisplayshortskip=7pt
	\section{Introduction}
	
	Monge-Amp\`{e}re equations, as a class of fully nonlinear equations, play important roles in many fields of analysis and geometry. 
	Many significant contributions have been made on their various aspects. In particular, for the Monge-Amp\`{e}re equation 
	\begin{equation}\label{eqma}
		\det(D^2v)=f\ \ \ \ \ \text{in}\ \mathbb{R}^n, 
	\end{equation}
	the asymptotic expansion of its solution near infinity has been widely discussed in the past few decades.  
	This can be regarded as an extension of Liouville's theorem.
	
	The classical theorem, in the case $f\equiv1$, states that any convex classical solution of \eqref{eqma} must be a quadratic 
	polynomial. This was obtained by J\"{o}rgens \cite{j}, Calabi \cite{ce} and Pogorelov \cite{p} for different $n$. Different 
	proofs and extensions can be found in \cite{cy,jx}. This theorem was extended to viscosity solution by Caffarelli \cite{ca}.
	
	If $f-1$ has a compact support, Caffarelli and Li \cite{cl} showed that any convex viscosity solution $v$ is very close to a 
	quadratic polynomial at infinity for $n\geq3$, while for $n=2$ there is an additional logarithmic term. In $\rr^2$, the 
	problem was studied by complex variable methods in \cite{d,fmm1,fmm2}. When $f$ is periodic, the relevant consequences were 
	established by Caffarelli-Li \cite{cl2} with positive $f\in C^\alpha$ and by Li-Lu \cite{ll} with positive $f\in L^\infty$. 
	The case $f$ is asymptotically periodic was considered by Teixeira-Zhang \cite{tz}. 
	
	{Comparing to the work in \cite{cl}, where $f$ is assumed to be $1$ outside a compact subset of $\rn$, 
		Bao, Li and Zhang \cite{blz} considered the case $f$ is a perturbation of $1$ near infinity.} 
	They established the similar asymptotic expansion under the assumptions on regularity and decay rate that 
	\begin{align}\label{eqc}
		&\exists m\geq3\ \text{such that $D^mf$ exists outside a compact subset of }\rn,\nonumber\\ 
		&\exists\beta>2\ \text{such that}\lim_{|x|\rightarrow+\infty}|x|^{\beta+k}|D^k(f(x)-1)|<+\infty,\ k=0,1,\cdots,m.
	\end{align}
	In addition, an example was given in \cite{blz} to 
	show that the decay rate $\beta>2$ is optimal.
	In the subsequent work \cite{lb}, Bao and Liu reduced $m$ from $3$ to $2$ by the detailed analysis of the solutions of 
	nonhomogeneous linearized equations. 
	
	A natural question is whether the regularity of $f$ can be reduced to $m<2$ and we will concentrate on it in this paper. In 
	fact, we just need that $m>0$ in our proof. To be more precise, our assumptions on $f$ are as follows:
	\vspace{-0.2cm}
	\begin{itemize}
		\item[(H)] The function $f\in C(\rn)\cap C_{loc}^{\alpha}({\rn\setminus\overline{\O}})$ for 
		some bounded open subset $\O\subset\rn$ and satisfies  
		\begin{equation}\label{eqfcdinfi}
			\limsup_{|x|\rightarrow+\infty}\left(|x|^{\beta}|f(x)-1|+
			|x|^{\beta+\az}[f]_{\text{{$C^\az(\overline{B_{|x|/2}(x)})$}}}\right)<+\infty
		\end{equation}
		for some $\alpha\in(0,1)$ and $\beta>2$.
	\end{itemize}
	\vspace{-0.2cm}
	Let $\mm^{n\times n}$ be the set of the real valued, $n\times n$ matrices and 
	\begin{equation*}
		\A:=\{A\in\mm^{n\times n}:\ A\text{ is symmetric, positive definite and }\det(A)=1\}.
	\end{equation*}
	Our main result is
	\begin{theorem}\label{thm1}
		Let $v\in C_{loc}(\rn)$ be a convex viscosity solution of \eqref{eqma}, where $f:\rn\rightarrow\rr$ satisfies 
		{\textnormal{(H)}}. If $n\geq3$, then there exist $c\in\rr$, $b\in\rn$ and $A\in\A$ such that
		\begin{equation*}
			\left\{
			\begin{aligned}
				&\limsup_{|x|\rightarrow+\infty}|x|^{\min\{n,\beta\}+k-2}
				\left|D^k\left(v(x)-\left(\frac{1}{2}x'Ax+b\cdot x+c\right)\right)\right|<+\infty,\ \ \ \beta\neq n,\\
				&\limsup_{|x|\rightarrow+\infty}|x|^{n+k-2}(\ln |x|)^{-1}
				\left|D^k\left(v(x)-\left(\frac{1}{2}x'Ax+b\cdot x+c\right)\right)\right|<+\infty,\ \ \ \beta=n.
			\end{aligned}
			\right.
		\end{equation*}
		Here, $k=0,1,2,$ and $x'$ denote the transpose of vector $x\in\rn$.
	\end{theorem}
	
	\begin{remark}\label{rem1}
		The decay \eqref{eqfcdinfi} in assumption \textnormal{(H)} is weaker than \eqref{eqc} even if $f$ is $C^1$. 
		In fact, it is easy to derive 
		\eqref{eqfcdinfi} from \eqref{eqc}. On the other hand, $Df_1(x)$ does not admit a limit at infinity for 
		$f_1(x):=1+e^{-|x|}\sin(e^{|x|})$, while \eqref{eqfcdinfi} holds for $f_1$. In fact, for any $|x|>1$ and 
		$z_1,z_2\in\overline{B_{|x|/2}(x)}$ (supposing that $|z_2|\leq|z_1|$), one has
		\begin{align*}
			\frac{|f_1(z_1)-f_1(z_2)|}{|z_1-z_2|^\alpha}&\leq|\sin(e^{{|z_2|}})|\frac{|e^{-|z_1|}-e^{-|z_2|}|}{|z_1-z_2|^\alpha}
			+e^{-|z_1|}\frac{|\sin(e^{{|z_1|}})-\sin(e^{{|z_2|}})|}{|z_1-z_2|^\alpha}:=J_1+J_2.
		\end{align*}
		It is easy to check that $J_1\leq Ce^{-|x|/2}|z_1-z_2|^{1-\alpha}$
		for constant $C$ independent of $x$. Note that $|\sin(e^{{|z_1|}})-\sin(e^{{|z_2|}})|\leq Ce^{|\xi|}|z_1-z_2|$ for some 
		$\xi=z_2+\theta(z_1-z_2)$ with $\theta\in(0,1)$ and $|\xi|\leq|z_1|$, we have 
		\begin{equation*}
			J_2\leq Ce^{-|z_1|}e^{|\xi|}|z_1-z_2|^{1-\alpha}\leq Ce^{-\frac{(1-\alpha)}{2}|x|}
			\ \ \ \ \text{if}\ |z_1-z_2|\leq e^{-\frac{|x|}{2}}
		\end{equation*}
		If $|z_1-z_2|\geq e^{-|x|/2}$, then $J_2\leq e^{-|z_1|}e^{\alpha|x|/2}\leq e^{-(1-\alpha)|x|/2}$.
		Thus $$[f_1]_{C^\alpha(\overline{B_{|x|/2}(x)})}\leq Ce^{-\frac{(1-\alpha)}{2}|x|}|x|^{1-\alpha},$$
		and $f_1$ verifies \eqref{eqfcdinfi} for any $\beta>2$.
	\end{remark}
	
	\begin{remark}\label{remf}
		It is clear that $f\in C^\alpha(\rn\setminus\O)$ under the assumption \textnormal{(H)}.
	\end{remark}
	
	\begin{remark}\label{remr}
		The viscosity solution $v$ in Theorem \ref{thm1} is really $C^{2,\alpha}$ near infinity and we will prove the regularity of 
		$v$ in Section 3.
	\end{remark}

	In this paper, a non-local method is developed to prove Theorem \ref{thm1}. 
	As mentioned above, the right term $f$ is at least $C^2$ in the 
	previous work and, of course, the derivative of $f$ plays a crucial role in their proof. Our method is essentially different 
	since the $f$ in this paper is just H\"{o}lder continuous and one can not expect the existence of its derivative. We introduce 
	some non-local objects and develop a series of non-local arguments to deal with the problems caused by the insufficient smoothness 
	of $f$.
	\subsection*{Notations and the structure of the article}
	
	Throughout the paper, we use the following notations:
	
	The point $x\in\rn$ will also be written as $x=(x^{(1)},\cdots,x^{(n)})$. 
	$B_r(x)\subset\rn$ denote the ball centered at $x$ with radius $r$. 
	We drop the center if it coincides with the origin, i.e., $B_r=B_r(0)$. 
	
	For a $n\times n$ matrix $B=(b_{ij})_{n\times n}$, $1\leq i,j\leq n$, $cof_{ij}B$ stands for the algebraic cofactor of 
	element $b_{ij}$. The identity matrix is denoted by $I$. 
	
	
	For multi-indices 
	$\gamma=(\gamma_1,\cdots,\gamma_n)\in\mathbb{N}^n$, we let $|\gamma|$ denote the sum of its components. 
	Given a function $u:\rn\rightarrow\rr$ and a point $x\in\rn$, we denote by $Du(x)$ and $\Delta u(x)$ the gradient  vector and 
	the  Laplacian respectively.  
	
	Let $k\in\mathbb{N}$, $\kappa\in(0,1)$ and $\Omega\subset\rn$ be open. 
	The space $C^k(\overline{\Omega})$ consists of functions $u: \overline{\Omega}\rightarrow\rr$, which admit derivatives 
	up to order $k$, such that 
	\begin{equation*}
		\|u\|_{C^{k}(\overline{\Omega})}:=\sup_{\substack{ \gamma\in\mathbb{N}^n\\|\gamma|\leq k}}
		\sup_{x\in\overline{\Omega}}\left|D^\gamma u(x)\right|<+\infty.
	\end{equation*}
	The H\"{o}lder space $C^{k,\kappa}(\overline{\Omega})$ consists of function $u\in C^k(\overline{\Omega})$ satisfying 
	\begin{gather*}
		\|u\|_{C^{k,\kappa}(\overline{\Omega})}:=\|u\|_{C^{k}(\overline{\Omega})}+[u]_{C^{k,\kappa}(\overline{\Omega})}<+\infty, 
	\end{gather*} 
	where the seminorm 
	\begin{equation*}
		[u]_{C^{k,\kappa}(\overline{\Omega})}:=\sup_{\substack{ \gamma\in\mathbb{N}^n\\|\gamma|=k}}
		\sup_{\substack{x,y\in\overline{\Omega} \\ x\neq y}}\frac{|D^\gamma u(x)-D^\gamma u(y)|}{|x-y|^\kappa}.
	\end{equation*}
	A function $u\in C^{k,\kappa}_{loc}({\Omega})$ if $u\in C^{k,\kappa}(K)$ for all compact 
	$K\subset\Omega$. 
	For convenient, we may write $C^{0}(\overline{\Omega})$, $C^{k,\kappa}(\overline{\Omega})$ as $C(\overline{\Omega})$, 
	$C^{k+\kappa}(\overline{\Omega})$,  respectively. 
	
	If $\det(D^2v(x))=f(x)$ in the classical sense for functions $v$ and $f$, then, by the property of determinant, we have 
	\begin{equation*}
		\sum_{1\leq i,j\leq n}cof_{ij}(D^2v(x))\partial_{ij}v(x)=nf(x).
	\end{equation*}
	Let $a_{ij}(x):=n^{-1}cof_{ij}(D^2v(x))$, then 
	\begin{equation}\label{notaeq}
		\sum_{1\leq i,j\leq n}a_{ij}(x)\partial_{ij}v(x)=f(x).
	\end{equation}
	For convenient, we will omit the notation $\sum$ and write \eqref{notaeq} as $a_{ij}(x)\partial_{ij}v(x)=f(x)$ in this article. 
	
	The paper is organized as follows: In Section 2, we introduce the fractional Laplacian and establish a series of crucial properties,  
	which will be applied in the sequel. In Section 3, we first show that the solution $v$ is locally $C^{2,\alpha}$ outside a bounded 
	set and then discuss the decay rate near infinity. Theorem \ref{thm1} is finally proved by a non-local method.

	\section{Fractional Laplacian and its properties}
	
	In this section we introduce the fractional Laplacian and some relevant objects. A series of crucial properties are then  
	established  for them. 
	
	
	The fractional Laplacian for a function $u:\rn\rightarrow\rr$ with the parameter $s\in(0,1)$ is defined as
	\begin{equation*}
		(-\Delta)^su(x)=c_{n,s} \textnormal{P.V.}\int_{\rn}\frac{u(x)-u(y)}{|x-y|^{n+2s}}dy,
	\end{equation*}
	where $c_{n,s}$ is the normalization constant and P.V. stands for the Cauchy principal value. 
	If $u$ is $C^{2s+\ep}$ at a point $x_0\in\rn$ and satisfies 
	\begin{equation}\label{eqfwdc}
		\|u\|_{L_s(\rn)}:=\int_{\rn}\frac{|u(y)|}{1+|y|^{n+2s}}dy<+\infty,
	\end{equation}
	then $(-\Delta)^su(x_0)$ can be calculated classically. 
	
	The fundamental solution of $(-\Delta)^s$ is 
	\begin{equation*}
		\Phi_s(x):=c_{n,-s}\frac{1}{|x|^{n-2s}},\ \ \ \ n>2s.
	\end{equation*}
	As in \cite{s}, we use the notation 
	\begin{equation*}
		(-\Delta)^{-s}u(x)=c_{n,-s}\int_{\rn}\frac{u(y)}{|x-y|^{n-2s}}dy, 
	\end{equation*}
	which is really the Riesz potential of $u$. If $u\in L^p(\rn)$ and $2sp>n$, then $(-\Delta)^{-s}u\in C(\rn)$, see \cite{ah,dpv}. 
	
	For reading convenient, we list some important properties for $(-\Delta)^s$ established in \cite[Section 2]{s}, which will 
	be frequently used in our paper. We first summarize \cite[Propositions 2.5-2.7]{s} as a following version that 
	can be directly applied.   
	\begin{proposition}\label{prof1}
		Let $s\in(0,1)$ and $u\in C^{k,\alpha}(\rn)$ for $\alpha\in(0,1)$ and $k=0,1,2$.
		\vspace{-0.2cm}
		\begin{itemize}
			\setlength{\itemsep}{-3pt}
			\item If $\alpha>2s$, then $(-\Delta)^su\in C^{k,\alpha-2s}(\rn)$ and 
			\begin{equation*}
				\left[(-\Delta)^su\right]_{C^{k,\alpha-2s}(\rn)}\leq C[u]_{C^{k,\alpha}(\rn)}.
			\end{equation*}
			\item If $\alpha<2s<1+\alpha$ and $k=1,2$, then $(-\Delta)^su\in C^{k-1,\alpha-2s+1}(\rn)$ and 
			\begin{equation*}
				\left[(-\Delta)^su\right]_{C^{k-1,\alpha-2s+1}(\rn)}\leq C[u]_{C^{k,\alpha}(\rn)}.
			\end{equation*}
		\end{itemize}
		\vspace{-0.2cm}
		The constant $C$ depends only on $\alpha,\ s$ and $n$.
	\end{proposition}
	
	The next property can be found in \cite[Proposition 2.8]{s}.
	
	\begin{proposition}\label{prof2}
		Let $s,\alpha\in(0,1)$ be such that $\alpha+2s<1$. If $u\in L^\infty(\rn)$ and $(-\Delta)^su\in C^{\alpha}(\rn)$, 
		then $u\in C^{\alpha+2s}(\rn)$ and 
		\begin{equation*}
			\|u\|_{C^{\alpha+2s}(\rn)}\leq C(\|u\|_{L^\infty(\rn)}+\|(-\Delta)^su\|_{C^{\alpha}(\rn)})
		\end{equation*}
		for a constant $C$ depending only on $\alpha,\ s$ and $n$.
	\end{proposition}
	
	More consequences about fractional Laplacian can be found in \cite{ah,b,clm,dpv,s}. 
	Now we begin to discuss the properties of fractional Laplacian. The following lemma is a direct calculation, see \cite[Page 432]{pp}.
	\begin{lemma}\label{lemmul}
		Let $\alpha\in(0,1)$ and $s\in(0,\alpha/2)$. If $u_1,u_2\in C^{\alpha}(\rn)$, then 
		\begin{align*}
			(-\Delta)^s(u_1u_2)(x)=&u_1(x)(-\Delta)^su_2(x)+u_2(x)(-\Delta)^su_1(x)\\
			&-c_{n,s}\int_{\rn}\frac{(u_1(x)-u_1(y))(u_2(x)-u_2(y))}{|x-y|^{n+2s}}dy.
		\end{align*}
	\end{lemma}
	
	The following result indicates that one can exchange the order of $(-\Delta)^s$ and $\partial$ when applying them to a suitable 
	function.
	\begin{lemma}\label{lemexc}
		Let $\alpha\in(0,1)$ and $s\in(0,1)$. If either of the following holds
		\vspace{-0.2cm}
		\begin{itemize}
			\item[(i)] $u\in C^{1,\alpha}(\rn)$ and $2s<\alpha$,
			\item[(ii)] $u\in C^{2,\alpha}(\rn)$ and $2s<1+\alpha$,
		\end{itemize}
		\vspace{-0.2cm}
		then for $k=1,2,\cdots,n$
		\begin{equation*}
			(-\Delta)^s(\partial_ku)(x)=\partial_k((-\Delta)^su)(x).
		\end{equation*}
	\end{lemma}
	
	\noindent{\it Proof.} In both case, we can apply Proposition \ref{prof1} to get that $(-\Delta)^su\in C^{1,\ep_1}(\rn)$ for 
	some $\ep_1>0$. On the other hand, $(-\Delta)^s(\partial_ku)(x)$ can be calculated classically and the integral converges 
	uniformly on $x$. Thus we can exchange the order of $(-\Delta)^s$ and $\partial$, as desired.
	\qed
	
	If $u\in C^1(\rn)$ satisfies
	\begin{equation*}\label{eqcfdc}
		\sup_{y\in B_{|x|/2}(x)}\frac{|u(x)-u(y)|}{|x-y|}\leq c|x|^{-\sigma-1},\ \ \ \text{for}\ |x|>R,
	\end{equation*}
	then, by taking $y\rightarrow x$, one immediately has $|Du(x)|\leq c|x|^{-\sigma-1}$ for $|x|>R$. 
	We next show a non-local analogue of this result.
	

	\begin{lemma}\label{lemd}
		Let $\alpha\in(0,1)$ and $s\in(0,\alpha/2)$. If $u\in C^{\alpha}(\rn)$ satisfies
		\begin{equation}\label{eqcfd2.5}
			|u(x)|\leq c'|x|^{-\sigma},\ \ [u]_{C^\alpha(\overline{{B_{|x|/2}(x)}})}\leq c'|x|^{-\sigma-\alpha},\ \ \ \text{for}\ |x|>R'
		\end{equation}
		for some positive constants $c',\sigma$ and $R'$, then there exists $C_1>0$ such that 
		\begin{equation*}
			|(-\Delta)^su(x)|\leq 
			\left\{
			\begin{aligned}
				C_1|x|^{-\min\{\sigma,n\}-2s},\ \ \ \sigma\neq n,\\
				C_1|x|^{-n-2s}(\ln|x|),\ \ \ \sigma=n,
			\end{aligned}
			\right.\ \ \ \ \ \ \text{for}\ |x|>2R',
		\end{equation*}    
		where the constant $C_1$ depends only on $n,\ s,\ \sigma,\ c'$, $R'$ and $\|u\|_{C(\rn)}$.
	\end{lemma}
	
	\noindent{\it Proof.} Without loss of generality, we assume $R'>1$. Define for $|x|>2R'$ that 
	\begin{align*}
		&A_1:=\{y\in\rn:\ |y|\leq|x|/2\},\\
		&A_2:=\{y\in\rn:\ |x-y|\leq|x|/2\},\\
		&A_3:=\rn\setminus\left(A_1\cup A_2\right).
	\end{align*}
	Then 
	\begin{equation*}
		(-\Delta)^su(x)
		=c_{n,s}\int_{A_1\cup A_2\cup A_3}\frac{u(x)-u(y)}{|x-y|^{n+2s}}dy.
	\end{equation*}
	
	Note that $|x-y|\geq|x|/2\geq|y|$ in $A_1$ and $|x|/2>1$, we write  
	\begin{align*}
		\int_{A_1}\frac{|u(x)-u(y)|}{|x-y|^{n+2s}}dy
		&\leq\frac{C}{|x|^{n+2s}}\int_{B_{|x|/2}}|u(x)-u(y)|dy\\
		&\leq\frac{C}{|x|^{n+2s}}\left(\int_{B_{R'}}+\int_{B_{|x|/2}\setminus B_{R'}}\right)|u(x)-u(y)|dy.
	\end{align*}
	Since $u$ is bounded, one can find a constant $C'>0$, which may depend on $\|u\|_{C(\rn)}$ and $R'$ such that   
	\begin{equation*}
		\int_{B_{R'}}|u(x)-u(y)|dy
		\leq C'.
	\end{equation*}
	On the other hand, equation \eqref{eqcfd2.5} yields that 
	\begin{align*}
		\int_{B_{|x|/2}\setminus B_{R'}}|u(x)-u(y)|dy\leq C\int_{B_{|x|/2}\setminus B_{R'}}|y|^{-\sigma}dy\leq
		\left\{
		\begin{array}{ll}
			C|x|^{n-\sigma},&\ \ \ \sigma\neq n,\\
			C\ln|x|,&\ \ \ \sigma=n.
		\end{array}
		\right.
	\end{align*}
	Thus 
	\begin{align*}
		\int_{B_{|x|/2}}|u(x)-u(y)|dy
		&\leq\left\{
		\begin{array}{ll}
			C(|x|^{n-\sigma}+1),&\ \ \ \sigma\neq n,\\
			C(\ln|x|+1),&\ \ \ \sigma=n.
		\end{array}
		\right.
	\end{align*}
	Note that 
	\begin{align*}
		&(|x|^{n-\sigma}+1)\cdot{|x|^{-n-2s}}\leq|x|^{-\min\{\sigma,n\}-2s},\\
		&(\ln|x|+1)\cdot{|x|^{-n-2s}}\leq|x|^{-n-2s}(\ln|x|),
	\end{align*}
	we obtain    
	\begin{align*}
		\int_{A_1}\frac{|u(x)-u(y)|}{|x-y|^{n+2s}}dy
		&\leq\left\{
		\begin{aligned}
			C|x|^{-\min\{\sigma,n\}-2s},\ \ \ \sigma\neq n,\\
			C|x|^{-n-2s}(\ln|x|),\ \ \ \sigma=n,
		\end{aligned}
		\right.\ \ \ \text{for}\ |x|>2R'.
	\end{align*}
	
	In $A_2$ we have $|x-y|\leq|x|/2\leq|y|$. Then \eqref{eqcfd2.5} allows us to calculate that  
	\begin{align*}
		\int_{A_2}\frac{|u(x)-u(y)|}{|x-y|^{n+2s}}dy
		\leq\frac{c'}{|x|^{\sigma+\alpha}}\int_{|x-y|\leq|x|/2}\frac{1}{|x-y|^{n+2s-\alpha}}dy
		\leq C|x|^{-\sigma-2s}.
	\end{align*}
	
	For $A_3$, we divide it into two parts
	\begin{equation*}
		A_3^+:=\{y\in A_3:\ |x-y|\geq|y|\}\ \ \ \text{and}\ \ \ A_3^-:=A_3\setminus A_3^+.
	\end{equation*}
	Then we have $|x|/2\leq|y|\leq|x-y|$ in $A_3^+$ and hence $|u(x)|+|u(y)|\leq C|x|^{-\sigma}$, which gives 
	\begin{align*}
		\int_{A_3^+}\frac{|u(x)-u(y)|}{|x-y|^{n+2s}}dy\leq
		\frac{C}{|x|^{\sigma}}\int_{|y|\geq|x|/2}\frac{1}{|y|^{n+2s}}dy
		\leq C|x|^{-\sigma-2s}.
	\end{align*}
	Finally, in $A_3^-$ there holds $|x|/2\leq|x-y|\leq|y|$ and thus we have 
	\begin{align*}
		\int_{A_3^-}\frac{|u(x)-u(y)|}{|x-y|^{n+2s}}dy\leq
		\frac{C}{|x|^{\sigma}}\int_{|x-y|\geq|x|/2}\frac{1}{|x-y|^{n+2s}}dy
		\leq C|x|^{-\sigma-2s}.
	\end{align*}
	The lemma follows from the above estimates.
	\qed
	
	\begin{remark}
		The result in Lemma \ref{lemd} is in some sense weaker than that for first order derivative. It is since $(-\Delta)^s$ 
		is non-local and the calculation of $(-\Delta)^su$ need the information of $u$ near the origin.
	\end{remark}
	

	Using the same argument as above we obtain the following corollary.
	\begin{corollary}\label{cor1}
		Let $\alpha_1,\alpha_2\in(0,1)$ and $s\in(0,(\alpha_1+\alpha_2)/2)$. If $u_i\in C^{\alpha_i}(\rn)$ satisfies 
		\begin{gather*}
			|u_i(x)|\leq c_i|x|^{-\sigma_i},\ \ [u_i]_{C^{\alpha_i}(\overline{B_{|x|/2}(x)})}\leq c_i|x|^{-\sigma_i-\alpha_i},
			\ \ \ \text{for }|x|>R',
		\end{gather*}
		with some positive constants $c_i,\sigma_i,R'$, $i=1,2$, then there exists $C>0$ such that 
		\begin{align*}
			&\left|\int_{\rn}\frac{(u_1(x)-u_1(y))(u_2(x)-u_2(y))}{|x-y|^{n+2s}}dy\right|\\
			&\leq 
			\left\{
			\begin{array}{ll}
				C|x|^{-\min\{\sigma_1+\sigma_2,n\}-2s},&\ \ \ \sigma_1+\sigma_2\neq n,\\
				C|x|^{-n-2s}(\ln|x|),&\ \ \ \sigma_1+\sigma_2=n,
			\end{array}
			\right.\ \ \ \ \ \ \text{for}\ |x|>2R',
		\end{align*} 
		where the constant $C$ depends only on $n,\ s,\ \sigma_i,\ c_i$, $R'$ and $\|u_i\|_{C(\rn)}$.   
	\end{corollary}
	
	
	In fact, if $u$ increases not too fast near infinity, then $(-\Delta)^su$ also admits a decay. 
	
	\begin{lemma}\label{lemi}
		Let $\alpha\in(0,1)$ and $s\in(0,\alpha/2)$. If $u$ is continuous and satisfies
		\begin{equation}\label{eqcfi}
			|u(x)|\leq c'|x|^{\sigma},\ \ [u]_{C^\alpha(\overline{B_{|x|/2}(x)})}\leq c'|x|^{\sigma-\alpha},\ \ \ \text{for}\ |x|>R'
		\end{equation}
		for some constants $c'>0,\ \sigma\in[0,2s)$ and $R'>0$, then there exist $C_1>0$ such that 
		\begin{equation*}
			|(-\Delta)^su(x)|\leq C_1|x|^{\sigma-2s}\ \ \ \ \ \text{for}\ |x|>2R',
		\end{equation*}    
		where the constant $C_1$ depends only on $n,\ s,\ \sigma,\ c'$, $R'$ and $\|u\|_{C(\overline{B_{R'}})}$.
	\end{lemma}
	
	\noindent{\it Proof.} The proof is similar to that of Lemma \ref{lemd}, we just list the differences.
	
	Let $A_1,A_2$ and $A_3$ be as in Lemma \ref{lemd}. In $A_1$, it is easy to check that 
	\begin{gather*}
		\int_{B_{R'}}|u(x)-u(y)|dy\leq C\leq	C|x|^{\sigma+n},\\
		\int_{B_{|x|/2}\setminus B_{R'}}|u(x)-u(y)|dy\leq C\int_{B_{|x|/2}\setminus B_{R'}}|y|^{\sigma}dy\leq	C|x|^{\sigma+n},
	\end{gather*}
	and thus 
	\begin{align*}
		\int_{A_1}\frac{|u(x)-u(y)|}{|x-y|^{n+2s}}dy\leq C|x|^{\sigma+n-n-2s}=C|x|^{\sigma-2s}\ \ \ \text{for}\ |x|>2R'.
	\end{align*}
	
	In $A_2$ we just need to use \eqref{eqcfi} to get 
	\begin{align*}
		\int_{A_2}\frac{|u(x)-u(y)|}{|x-y|^{n+2s}}dy
		\leq c'|x|^{\sigma-\alpha}\int_{|x-y|\leq|x|/2}\frac{1}{|x-y|^{n+2s-\alpha}}dy
		\leq C|x|^{\sigma-2s}.
	\end{align*}
	
	Since $\sigma<2s$, combine with the fact $|u(x)-u(y)|\leq C|y|^{\sigma}$ one has 
	\begin{align*}
		\int_{A_3^+}\frac{|u(x)-u(y)|}{|x-y|^{n+2s}}dy\leq
		C\int_{|y|\geq|x|/2}\frac{|y|^{\sigma}}{|y|^{n+2s}}dy
		\leq C|x|^{\sigma-2s}.
	\end{align*}
	Finally, in $A_3^-$ there holds $|x|/2\leq|x-y|\leq|y|$. By the triangle inequality we also have 
	\begin{equation*}
		|y|\leq|x-y|+|x|\leq3|x-y|.
	\end{equation*}
	Thus, 
	\begin{align*}
		\int_{A_3^-}\frac{|u(x)-u(y)|}{|x-y|^{n+2s}}dy&\leq
		C\int_{|x-y|\geq|x|/2}\frac{|y|^{\sigma}}{|x-y|^{n+2s}}dy\\
		&\leq C\int_{|x-y|\geq|x|/2}\frac{|x-y|^{\sigma}}{|x-y|^{n+2s}}dy
		\leq C|x|^{\sigma-2s}.
	\end{align*}
	The proof is completed.
	\qed

	Next we establish the decay for Riesz potential.
	
	\begin{corollary}\label{cor2}
		Let $s\in(0,1)$ and $u\in C(\rn)$ be such that 
		\begin{equation*}
			|u(x)|\leq c'|x|^{-\sigma},\ \ \ \text{for}\ |x|>R'
		\end{equation*}
		for some positive constants $c',\sigma>2s$ and $R'$. Then 
		\begin{equation*}
			|(-\Delta)^{-s}u(x)|\leq 
			\left\{
			\begin{aligned}
				C|x|^{2s-\min\{\sigma,n\}},\ \ \ \sigma\neq n,\\
				C|x|^{2s-n}(\ln|x|),\ \ \ \sigma=n,
			\end{aligned}
			\right.\ \ \ \ \ \ \text{for}\ |x|>2R',
		\end{equation*} 
		where the constant $C>0$ depends only on $n,\ s,\ \sigma,\ c'$, $R'$ and $\|u\|_{C(\rn)}$.   
	\end{corollary}
	
	\noindent{\it Proof.} It is easy to check that $u\in L^p(\rn)$ for $p$ large enough, 
	thus $(-\Delta)^{-s}u\in C(\rn)$, see \cite[Theorem 1.2.4]{ah}. 
	Let $A_1,\ A_2$ and $A_3$ be as in Lemma \ref{lemd}, then 
	\begin{equation*}
		(-\Delta)^{-s}u(x)=c_{n,-s}\int_{\rn}\frac{u(y)}{|x-y|^{n-2s}}dy
		=c_{n,-s}\int_{A_1\cup A_2\cup A_3}\frac{u(y)}{|x-y|^{n-2s}}dy.
	\end{equation*}
	The estimate in $A_1$ is similar to that of Lemma \ref{lemd}, one just needs to replace $s$ by $-s$. 
	In $A_2$ we have $|y|\geq|x|/2\geq|x-y|$ and thus 
	\begin{align*}
		\int_{A_2}\frac{|u(y)|}{|x-y|^{n-2s}}dy
		&\leq\frac{C}{|x|^{\sigma}}\int_{B_{|x|/2}(x)}\frac{1}{|x-y|^{n-2s}}dy
		\leq C|x|^{2s}\cdot|x|^{-\sigma}\leq C|x|^{2s-\sigma}.
	\end{align*}
	To compute the integral in $A_3$, we also use the notations $A_3^+$ and $A_3^-$ as in Lemma \ref{lemd}. Note that we have   
	$|u(y)|\leq c'|y|^{-\sigma}$ and $|x-y|\geq|y|\geq|x|/2$ in $A_3^+$, which yield  
	\begin{align*}
		\int_{A_3^+}\frac{|u(y)|}{|x-y|^{n-2s}}dy\leq
		C\int_{|y|\geq|x|/2}\frac{1}{|y|^{n-2s+\sigma}}dy
		\leq C|x|^{2s-\sigma}
	\end{align*}
	since $\sigma>2s$. In $A_3^-$ the inequalities $|x|/2\leq|x-y|\leq|y|$ hold, thus 
	$|u(y)|\leq c'|y|^{-\sigma}\leq c'|x-y|^{-\sigma}$ and 
	\begin{align*}
		\int_{A_3^-}\frac{|u(y)|}{|x-y|^{n+2s}}dy\leq
		C\int_{|x-y|\geq|x|/2}\frac{1}{|x-y|^{n-2s+\sigma}}dy
		\leq C|x|^{2s-\sigma}.
	\end{align*}
	The proof is completed. 
	\qed
	
	We end up this section with the following H\"{o}lder continuity result.
	
	\begin{lemma}\label{lemeca}
		Let $\alpha_i\in(0,1)$, $u_i\in C^{\alpha_i}(\rn)$, $i=1,2$, and $s\in(0,\min\{\alpha_1,\alpha_2\}/2)$. Then 
		\begin{equation*}
			I(x):=\int_{\rn}\frac{(u_1(x)-u_1(y))(u_2(x)-u_2(y))}{|x-y|^{n+2s}}dy\in C^{\min\{\alpha_1,\alpha_2\}-2s}(\rn)
		\end{equation*}   
		and 
		\begin{equation*}
			\|I\|_{C^{\min\{\alpha_1,\alpha_2\}-2s}(\rn)}\leq C_1\|u_1\|_{C^{\alpha_1}(\rn)}\|u_2\|_{C^{\alpha_2}(\rn)}
		\end{equation*}
		for constant $C_1$ depending only on $n$ and $s$.
	\end{lemma}
	\noindent{\it Proof.} Under the assumptions on $u_i$, $I(x)$ is clearly bounded in $\rn$. In fact, one can compute that 
	\begin{align*}
		|I(x)|\leq&\left(\int_{B_1(x)}+\int_{\rn\setminus B_1(x)}\right)\frac{(u_1(x)-u_1(y))(u_2(x)-u_2(y))}{|x-y|^{n+2s}}dy\\
		\leq&\int_{B_1(x)}\frac{2[u_1]_{C^{\alpha_1}(\rn)}\|u_2\|_{C(\rn)}}{|x-y|^{n+2s-\alpha_1}}dy
		+\int_{\rn\setminus B_1(x)}\frac{4\|u_1\|_{C(\rn)}\|u_2\|_{C(\rn)}}{|x-y|^{n+2s}}dy\leq C.
	\end{align*}
	
	Without loss of generality, we assume that $\alpha_1\leq\alpha_2$. Then it is enough to show that 
	\begin{equation}\label{eqi}
		\left[I(x)\right]_{C^{\alpha_1-2s}(\rn)}\leq C.
	\end{equation}
	By changing variables we have 
	\begin{equation*}
		I(x)= \int_{\rn}\frac{(u_1(x)-u_1(x-y))(u_2(x)-u_2(x-y))}{|y|^{n+2s}}dy.
	\end{equation*}
	Let $x_1,x_2\in\rn$, $r:=|x_1-x_2|$ and $g(x,y):=(u_1(x)-u_1(x-y))(u_2(x)-u_2(x-y))$, then 
	\begin{equation*}
		|I(x_1)-I(x_2)|\leq\left(\int_{|y|\leq r}+\int_{|y|>r}\right)\frac{|g(x_1,y)-g(x_2,y)|}{|y|^{n+2s}}dy:=D_1+D_2.
	\end{equation*}
	From the inequalities 
	\begin{equation*}
		|u_1(x)-u_1(x-y)|\leq [u_1]_{C^{\alpha_1}(\rn)}|y|^{\alpha_1}\ \ \ \text{and}\ \ \ 
		|u_2(x)-u_2(x-y)|\leq2\|u_2\|_{C(\rn)}
	\end{equation*}
	we obtain $|g(x,y)|\leq C'|y|^{\alpha_1}$ for $y\in\rn$, where the constant $C'$ depends on the $C^{\alpha_i}$-norm 
	of $u_i$ but not on $x$. Since $\alpha_1-2s\in(0,1)$, it follows 
	\begin{equation}\label{eqd1}
		D_1\leq\int_{|y|\leq r}\frac{C|y|^{\alpha_1}}{|y|^{n+2s}}dy\leq Cr^{\alpha_1-2s}=C|x_1-x_2|^{\alpha_1-2s}.
	\end{equation}
	To estimate $D_2$, we first claim that 
	\begin{equation}\label{eqclaim}
		[u_1(\cdot-z_1)u_2(\cdot-z_2)]_{C^{\alpha_1}(\rn)}\leq A
		\ \ \ \text{for}\ z_1,z_2\in\rn,
	\end{equation}
	where the constant $A>0$ is independent of $z_1$ and $z_2$. 
	This implies by a direct calculation that $|g(x_1,y)-g(x_2,y)|\leq C|x_1-x_2|^{\alpha_1}$. Thus, 
	\begin{equation}\label{eqd2}
		D_2\leq\int_{|y|>r}\frac{C|x_1-x_2|^{\alpha_1}}{|y|^{n+2s}}dy\leq Cr^{-2s}|x_1-x_2|^{\alpha_1}=C|x_1-x_2|^{\alpha_1-2s}.
	\end{equation}
	
	Now we turn to prove \eqref{eqclaim}. Let $x,y\in\rn$. If $|x-y|\leq1$, then 
	by inserting the term $u_1(y-z_1)u_2(x-z_2)$ we have 
	\begin{align*}
		&\frac{|u_1(x-z_1)u_2(x-z_2)-u_1(y-z_1)u_2(y-z_2)|}{|x-y|^{\alpha_1}}\\
		&\leq[u_1]_{C^{\alpha_1}(\rn)}\|u_2\|_{C(\rn)}+\|u_1\|_{C(\rn)}[u_2]_{C^{\alpha_2}(\rn)}\leq A.
	\end{align*}
	If $|x-y|\geq1$, then 
	\begin{equation*}
		\frac{|u_1(x-z_1)u_2(x-z_2)-u_1(y-z_1)u_2(y-z_2)|}{|x-y|^{\alpha_1}}\leq
		2\|u_1\|_{C(\rn)}\|u_2\|_{C(\rn)}\leq A.
	\end{equation*}
	Thus \eqref{eqclaim} is proved.
	
	Note that the constant $C$ in \eqref{eqd1} and \eqref{eqd2} can be write as 
	$C_1\|u_1\|_{C^{\alpha_1}(\rn)}\|u_2\|_{C^{\alpha_2}(\rn)}$ for constant $C_1$ depending only on $n$ and $s$. 
	We obtain \eqref{eqi} by adding \eqref{eqd1} and \eqref{eqd2}. The proof is completed.
	\qed

	\section{The proof of Theorem \ref{thm1}}
	
	The goal of this section is to prove Theorem \ref{thm1}.

	We first state a priori consequence without proof, which was established in {\cite[Theorem 1.2]{blz}} and \cite[Theorem 2.2]{lb}, 
	based on the level set method by [7, Section 3].
	
	\begin{proposition}\label{pro1}
		Under the assumptions in Theorem \ref{thm1}, there exists a linear transform $T$ satisfying $\det T=1$ such that
		\begin{equation*}
			\left|\left(v\circ T\right)(x)-\frac{1}{2}|x|^2\right|\leq c_1|x|^{2-\ep},\ \ \ \text{for}\ |x|>R_0
		\end{equation*}
		for some $c_1>0$, $R_0>1$ and $\ep>0$.
	\end{proposition}
	
	According to Proposition \ref{pro1}, we just need to establish Theorem \ref{thm1} under the assumption on $v$ that 
	there exist $c_1>0$, $\ep>0$ and $R_0>1$ such that
	\begin{equation}\label{eqv1}
		\left|v(x)-\frac{1}{2}|x|^2\right|\leq c_1|x|^{2-\ep},\ \ \ \ \text{for}\ |x|>R_0.
	\end{equation}
	For convenient, we use from now the notation 
	\begin{equation}\label{maeqw}
		w(x):=v(x)-\frac{1}{2}|x|^2.
	\end{equation}
	
	We will prove that $v$, and hence $w$, is $C_{loc}^{2,\alpha}$ near infinity and then establish 
	the decay of $D^kw$ for $k=0,1,2$. This is Lemma \ref{lem21}. We remind the reader that the interior estimates for the 
	solution to Monge-Amp\`{e}re equation can be found in \textnormal{\cite{c,fjm,jw}} and some of these results will be 
	applied in the following lemma.
	
	
	\begin{lemma}\label{lem21}
		Let $v$ and $f$ be as in Theorem \ref{thm1}. If $v$ verifies \eqref{eqv1}, then 
		there exist $C_1>0$ and $R_1>R_0$ such that $v\in C^{2,\az}_{loc}(\rn\setminus\overline{B_{R_1}})$ and 
		\begin{equation*}
			\left\{
			\begin{aligned}
				&|D^kw(x)|\leq C_1|x|^{2-k-\beta_{\ep}},\ \ \ k=0,1,2,\ \ |x|>R_1,\\
				&\frac{|D^2w(x_1)-D^2w(x_2)|}{|x_1-x_2|^{\alpha}}\leq C_1|x_1|^{-\beta_\ep-\alpha},\ \ \ |x_1|>R_1,
				\ |x_2-x_1|<\frac{|x_1|}{2},
			\end{aligned}
			\right.
		\end{equation*}
		where $\beta_\ep=\min\{\ep,\beta\}$. The constants $C_1$ and $R_1$ depend only on 
		$n,\ R_0,\ \ep,\ \az,\ \beta,\ c_1$ and $\|f\|_{C(\rn)}$.
	\end{lemma}
	
	\noindent{\it Proof.} For $|x|=R>2R_0$, we define
	\begin{gather*}
		v_R(y):=\left(\frac{4}{R}\right)^2v\left(x+\frac{R}{4}y\right),\ \ \ |y|\leq2,\\
		w_R(y):=\left(\frac{4}{R}\right)^2w\left(x+\frac{R}{4}y\right),\ \ \ |y|\leq2.
	\end{gather*}
	Then \eqref{eqv1} yields that 
	\begin{equation*}
		\|v_R\|_{C(\overline{B_2})}\leq C,\ \ \ \|w_R\|_{C(\overline{B_2})}\leq CR^{-\ep}.
	\end{equation*}
	Therefore, we have 
	\begin{equation*}
		v_R(y)-\left(\frac{1}{2}|y|^2+\frac{4}{R}x\cdot y+\frac{8}{R^2}|x|^2\right)=O(R^{-\ep}),\ \ \ |y|\leq2.
	\end{equation*}
	If we define 
	\begin{equation*}
		\bar{v}_R(y):=v_R(y)-\frac{4}{R}x\cdot y-\frac{8}{R^2}|x|^2,
	\end{equation*}
	then 
	the $\bar{v}_R$ satisfies in viscosity sense that
	\begin{equation}\label{eqvr1}
		\det(D^2\bar{v}_R(y))=f_R(y):=f\left(x+\frac{R}{4}y\right),\ \ \ |y|<2.
	\end{equation}
	We get from (H) that
	\begin{equation}\label{eqfr1}
		\begin{aligned}
			\|f_R-1\|_{C^\alpha(\overline{B_2})}\leq CR^{-\beta}
		\end{aligned}
	\end{equation}
	It is easy to check that 
	\begin{equation*}
		B_{7/5}\subset\Omega_1:=\{|y|<2:\ \bar{v}_R(y)<1\}\subset B_2
	\end{equation*}
	for $R>R_1$ with $R_1$ large enough. Moreover, $\Omega_1$ is convex and $\bar{v}_R-1=0$ on $\partial\Omega_1$. 
	Thus we can apply the interior estimate \cite[Theorem 2]{c} to $\bar{v}_R-1$ on the domain $\Omega_1$ to obtain 
	that $\bar{v}_R-1\in C^{2,\alpha}(\overline{B_{13/10}})$, which allows us to get from \cite[Corollary 4]{ca1} the 
	strict convexity of $\bar{v}_R-1$ 
	on $B_{13/10}$. Then the Schauder estimate \cite[Theorem 1.1]{fjm} gives
	\begin{equation*}
		\|D^2v_R\|_{C^\alpha(\overline{B_{6/5}})}=\|D^2\bar{v}_R\|_{C^\alpha(\overline{B_{6/5}})}\leq C.
	\end{equation*}
	We rewrite \eqref{eqvr1} as 
	\begin{equation*}
		a_{ij}^R\partial_{ij}v_R=f_R\ \ \ \ \text{in}\ B_2
	\end{equation*}
	with $a_{ij}^R=n^{-1}cof_{ij}(D^2v_R)$. Clearly, \eqref{eqvr1} and \eqref{eqfr1} yields that 
	\begin{equation}\label{eqvrc2b}
		\frac{I}{C}\leq D^2v_R\leq CI\ \ \ \ \text{in}\ B_{6/5}
	\end{equation}
	for some $C$ independent of $R$. Therefore, $a_{ij}^R$ is uniformly elliptic and 
	\begin{equation*}
		\|a_{ij}^R\|_{C^{\alpha}(\overline{B_{6/5}})}\leq C.
	\end{equation*}
	By Schauder estimate we have 
	\begin{equation*}
		\|v_R\|_{C^{2,\alpha}(\overline{B_1})}\leq C(\|v_R\|_{C(\overline{B_2})}+\|f_R\|_{C^{\alpha}(\overline{B_{2}})})\leq C.
	\end{equation*}
	
	Now we turn to estimate $w_R$. The difference between \eqref{eqvr1} and $\det(I)=1$ gives 
	\begin{equation*}
		\tilde{a}^R_{ij}\partial_{ij}w_R=f_R-1=O(R^{-\beta}),
	\end{equation*}
	where $\tilde{a}^R_{ij}(y)=\int_0^1cof_{ij}(I+tD^2w_R(y))dt$. By \eqref{eqvrc2b} we have 
	\begin{equation*}
		\frac{I}{C}\leq \left(\tilde{a}^R_{ij}\right)_{n\times n}\leq CI\ \ \ \ \text{on}\ B_{6/5}.
	\end{equation*}
	The Schauder estimate gives 
	\begin{equation*}\label{eqwrc2}
		\|w_R\|_{C^{2,\alpha}(\overline{B_1})}\leq 
		C(\|w_R\|_{C(\overline{B_2})}+\|f_R-1\|_{C^{\alpha}(\overline{B_2})})\leq CR^{-\beta_\ep}.
	\end{equation*}
	This completes the proof.
	\qed

	
	Since $v\in C^{2,\az}_{loc}(\rn\setminus\overline{B_{R_1}})$, we can use \cite[Theorem 3.4]{y} to redefine the value of $v$ 
	inside $B_{2R_1}$ such that the new function $\bar{v}$ is locally $C^{2,\az}$ on $\rn$. Let $\bar{f}:=\det(D^2\bar{v})$, then 
	$\bar{f}\in C^\alpha(\rn)$ coincides with $f$ outside $B_{2R_1}$. This allows us to assume that 
	$v\in C^{2,\az}_{loc}(\rn)$ and ${f}\in C^\alpha(\rn)$ from now, just to apply the fractional Laplacian conveniently. 
	
	We remind the reader that the redefinition of $v$ and thus $f$ may influences the constants $C_1$ and $R_1$ but has no effect 
	on $\A$, $b$ and $c$ in Theorem \ref{thm1}.
	
	It follows from $v\in C^{2,\az}_{loc}(\rn)$ that $w\in C^{2,\alpha}_{loc}(\rn)$. Combining with the decay on $D^2w$ established 
	in Lemma \ref{lem21}, one implies that $\partial_{ij}w\in C^{\az}(\rn)$, $1\leq i,j\leq n$.
	Note that $(-\Delta)^sw$ and $(-\Delta)^s\partial w$ are not always well defined for sufficiently small $s>0$ since, 
	by Lemma \ref{lem21}, $|w|$ and $|Dw|$ may increase faster than $|x|^{2s}$ near infinity and thus may 
	not satisfy \eqref{eqfwdc}.  
	
	We then give a useful result.
	
	\begin{lemma}\label{lemclaim}
		Let $u\in C^1(\rn)$ and 
		\begin{equation}\label{eqdclaim}
			|Du(x)|\leq M|x|^{-1-\sigma}\ \ \ \ \ \text{for}\ |x|\geq R',
		\end{equation} 
		with positive constants $M,\ \sigma$ and $R'$. Then, for some $u_0\in\rr$, 
		we have $|u(x)-u_0|\leq C|x|^{-\sigma}$ for $|x|>R'$, where the constant $C$ depends only on $M$ and $n$.
	\end{lemma}
	\noindent{Proof.} We first show by contradiction that $u$ is bounded at least on one side. 
	Indeed, if $u$ is unbounded on both, i.e., 
	\begin{equation*}
		\limsup_{|x|\rightarrow+\infty}u(x)=+\infty\ \ \ \text{and}\ \ \ \liminf_{|x|\rightarrow+\infty}u(x)=-\infty,
	\end{equation*}
	then, by the continuity of $u$, there exists a sequence $\{x_i\}_{i=1}^{+\infty}$ such that 
	$1<|x_i|<|x_{i+1}|\rightarrow+\infty$ and $u(x_i)=0$. We claim that
	\begin{equation}\label{eqdcb1}
		|u(x_0)|\leq C|x_i|^{-\sigma}
	\end{equation}
	for any $x_0\in\partial B_{|x_i|}$, where the constant $C$ depends only on $M$ and $n$.
	
	To show \eqref{eqdcb1}, we denote by $\Gamma$ the minor arc connecting $x_0$ and $x_i$ in the great circle of sphere 
	$\partial B_{|x_i|}$ and consider the parametric equation of $\Gamma$: 
	$x=x(l)$, where $l\in[0,L]$ is the arc-length parameter (the arc-length between $x_i$ and a 
	point $x\in\Gamma$) and $L$ is the length of $\Gamma$. Let $T(l)=u(x(l))$, then $T(0)=u(x_i)=0$, $T(L)=u(x_0)$ and 
	\begin{equation*}
		T(L)=T(0)+\int_0^L\frac{dT}{dl}dl=\int_0^L\frac{dT}{dl}dl.
	\end{equation*}
	To establish \eqref{eqdcb1}, we just need to show $|dT/dl|\leq C|x_i|^{-1-\sigma}$ since $L\leq2\pi|x_i|$. Observe
	\begin{equation*}
		\frac{dT}{dl}=Du\cdot\frac{dx}{dl},\ \ \ \left(\frac{dx}{dl}=(x^{(1)}_l,\cdots,x^{(n)}_l)\right),
	\end{equation*}
	by \eqref{eqdclaim}, it is enough to prove $|dx/dl|\leq C$. By rotating the axis (the new coordinate is denoted by $z$) 
	we can put the minor arc $\Gamma$ into the coordinate plane $z^{(1)}Oz^{(2)}$ with $x_i$ becomes the point $(|x_i|,0,\cdots,0)$. 
	The new equations of $\Gamma$ are 
	\begin{equation*}
		z^{(1)}=|x_i|\cos\frac{l}{|x_i|},\ z^{(2)}=|x_i|\sin\frac{l}{|x_i|},\ z^{(m)}=0\ (m\geq3),\ l\in[0,L].
	\end{equation*}
	It is easy to calculate that $|dz/dl|^2=\sum_{m=1}^n \{(z^{(m)})_l\}^2=1$. Since the rotation, which is a orthogonal 
	transformation, depends only on $x_i$ and $x_0$, we obtain $|dx/dl|\leq C_n$. 
	
	For any $x'\in B_{|x_{i+1}|}\setminus\overline{B_{|x_i|}}$, let $x'_i:=x'|x_i|/|x'|$. Then $x'_i\in\partial B_{|x_i|}$ and 
	\begin{equation}\label{eqdcb2}
		|u(x')-u(x'_i)|\leq C\left(|x_i|^{-\sigma}-|x_{i+1}|^{-\sigma}\right),
	\end{equation}
	which follows from a similar but easier process as above. 
	In fact, one just need to take $\Gamma$ as the line connecting $x'_i$ and $x'$:
	\begin{equation*}
		x(l)=\frac{x'}{|x'|}l,\ \ \ \ l\in\left[|x'_i|,|x'|\right].
	\end{equation*}
	Clearly, $|dx/dl|=1$.
	
	The inequalities \eqref{eqdcb1} and \eqref{eqdcb2} yield that $|u(x)|\leq C$ for $|x|>R'$ since $i$ is arbitrary, a 
	contradiction. Thus $u$ is bounded from at least one side. Without loss of generality, we assume that $u\geq M_0$. 
	
	Let now
	\begin{equation*}
		\liminf_{|x|\rightarrow+\infty}u(x)=u_0\geq M_0.
	\end{equation*}
	Then for any $\mu>0$ there exist a sequence $\{y_i\}_{i=1}^{+\infty}$ with $|y_i|\rightarrow+\infty$ such that 
	$|u(y_i)-u_0|<\mu$. The same argument as \eqref{eqdcb1} and \eqref{eqdcb2} shows that
	\begin{equation*}
		u(y)\leq u_0+\mu+C|y_i|^{-\sigma}\ \ \ \text{and}\ \ \ u(y')\leq u_0+\mu+2C|y_i|^{-\sigma}
	\end{equation*}
	for $y\in\partial B_{|y_i|}$ and $y'\in B_{|y_{i+1}|}\setminus\overline{B_{|y_i|}}$. Letting $y'$, and hence $y_i$, tend to 
	$+\infty$ we obtain 
	\begin{equation*}
		\limsup_{|y'|\rightarrow+\infty}u(y')\leq u_0+\mu.
	\end{equation*}
	By taking $\mu\rightarrow0^+$ we get $u(x)\rightarrow u_0$ as $|x|\rightarrow0$. 
	
	For any $|x|,|y|>R'$ we obtain by \eqref{eqdcb1} and \eqref{eqdcb2} that $|u(x)-u(y)|\leq C(|x|^{-\sigma}-|y|^{-\sigma})$.  
	Passing to the limit $y\rightarrow+\infty$ we get $|u(x)-u_0|\leq C|x|^{-\sigma}$.
	\qed
	
	\begin{remark}\label{remdec}
		If $u$ is continuous and replacing \eqref{eqdclaim} in Lemma \ref{lemclaim} by 
		\begin{equation*}
			|Du(x)|\leq M|x|^{-1+\sigma}\ \ \ \ \ \text{for}\ |x|\geq R',
		\end{equation*} 
		with positive constants $M,\ \sigma$ and $R'$. Then, similar to the argument in proving \eqref{eqdcb1} and \eqref{eqdcb2}, we 
		have $|u(x)|\leq C|x|^{\sigma}$ for $|x|\geq 1$.
	\end{remark}
	
	%
	
	Next we establish a faster decay than that in Lemma \ref{lem21} by applying the non-local operator $(-\Delta)^s$.
	
	\begin{lemma}\label{lem22}
		Under the assumption of Lemma \ref{lem21} and let $R_1$ be the large constant determined in the proof of Lemma \ref{lem21}.
		If in addition $2\ep<1$, then
		\begin{equation*}
			\left\{
			\begin{aligned}
				&|D^kw(x)|\leq C|x|^{2-k-2\ep},\ \ \ k=0,1,2,\ \ |x|>2R_1,\\
				&\frac{|D^2w(x_1)-D^2w(x_2)|}{|x_1-x_2|^{\alpha}}\leq C|x_1|^{-2\ep-\alpha},\ \ \ |x_1|>2R_1,
				\ |x_2-x_1|<\frac{|x_1|}{2}.
			\end{aligned}
			\right.
		\end{equation*}
	\end{lemma}
	
	\noindent{\it Proof.} 
	Note that the equation for $w$ is $\det(D^2w+I)=f$. As previous, $w$ satisfies 
	\begin{equation}\label{eq23v}
		\tilde{a}_{ij}\partial_{ij}w=f-1\ \ \ \ \text{in}\ \rn,
	\end{equation}
	where $\tilde{a}_{ij}(y)=\int_0^1cof_{ij}(I+tD^2w(y))dt$.  
	Then by Lemma \ref{lem21} we can compute that 
	\begin{equation*}
		|\tilde{a}_{ij}(x)-\delta_{ij}|\leq C|x|^{-\ep},\ \ \ |x|>R_1
	\end{equation*}
	and
	\begin{equation*}
		\frac{|\tilde{a}_{ij}(x_1)-\tilde{a}_{ij}(x_2)|}{|x_1-x_2|^\alpha}\leq C|x_1|^{-\ep-\alpha},\ \ |x_1|>R_1,\ \ |x_2-x_1|<|x_1|/2.
	\end{equation*}
	
	
	Let $s\in(0,\alpha/2)$. Note that $\tilde{a}_{ij},\ \partial_{ij}w$ and $f$ are all belong to $C^\az(\rn)$, by Lemma \ref{lemmul} 
	we can apply $(-\Delta)^s$ to \eqref{eq23v} to obtain that
	\begin{equation*}
		\begin{aligned}
			\tilde{a}_{ij}(-\Delta)^s\partial_{ij}w=&(-\Delta)^sf-\partial_{ij}w(-\Delta)^s\tilde{a}_{ij}\\
			&+c_{n,s}\int_{\rn}\frac{\left(\tilde{a}_{ij}(x)-\tilde{a}_{ij}(y)\right)\left(\partial_{ij}w(x)-\partial_{ij}w(y)\right)}
			{|x-y|^{n+2s}}dy.
		\end{aligned}
	\end{equation*}
	Then we write the above equation as 
	\begin{equation*}
		(-\Delta)^s(\Delta w)=F
	\end{equation*}
	where 
	\begin{equation*}
		\begin{aligned}
			F(x):=&(-\Delta)^sf(x)-\partial_{ij}w(x)(-\Delta)^s\tilde{a}_{ij}(x)
			-(\tilde{a}_{ij}(x)-\delta_{ij})(-\Delta)^s(\partial_{ij}w)(x)\\
			&+c_{n,s}\int_{\rn}\frac{\left(\tilde{a}_{ij}(x)-\tilde{a}_{ij}(y)\right)\left(\partial_{ij}w(x)-\partial_{ij}w(y)\right)}
			{|x-y|^{n+2s}}dy.
		\end{aligned}
	\end{equation*}
	Combined with Lemma \ref{lemd} and the assumptions on $f$ we obtain that
	\begin{equation*}
		|(-\Delta)^sf|\leq 
		\left\{
		\begin{aligned}
			C_1|x|^{-2s-\min\{\beta,n\}},\ \ \ \beta\neq n,\\
			C_1|x|^{-2s-n}(\ln|x|),\ \ \ \beta=n.
		\end{aligned}
		\right.
	\end{equation*}
	On the other hand, Lemma \ref{lem21} and Corollary \ref{cor1} imply that the other three terms in the definition of $F$ have 
	the same decay with exponent $-2s-2\ep$.
	Since $2\ep<1$  and $\beta>2$, it follows
	\begin{equation}\label{eqF}
		|F(x)|\leq C|x|^{-2s-2\ep},\ \ \ |x|>2R_1.
	\end{equation}
	By applying Lemma \ref{lemeca} and Proposition \ref{prof1} we also have $F\in C^{\az-2s}(\rn)$. 
	
	As introduced in Section 2, we consider the Riesz potential of $F$
	\begin{equation*}
		H(x):=c_{n,-s}\int_{\rn}\frac{F(y)}{|x-y|^{n-2s}}dy.
	\end{equation*}
	Note that $F$ is continuous, by \eqref{eqF} we see that $F\in L^p(\rn)$ for $p$ large enough such that $2sp>n$. Then 
	\cite[Theorem 1.2.4]{ah} implies $H\in C(\rn)$. We also obtain from \cite{s} that $(-\Delta)^{s}H(x)=F(x)$. By 
	Corollary \ref{cor2}, there holds  
	\begin{equation}\label{eqHl}
		|H(x)|\leq C|x|^{-2\ep},\ \ \ |x|>2R_1.
	\end{equation}
	Moreover, Proposition \ref{prof2} yields that $H\in C^{\alpha}(\rn)$. 
	
	Now we have 
	\begin{equation*}
		(-\Delta)^s\left(\Delta w-H\right)=0\ \ \ \ \text{in}\ \rn\setminus B_{2R_1}.
	\end{equation*}
	It is obvious that $\Delta w-H\rightarrow0$ at infinity. By comparing with a multiple of $|x|^{2s-n}$ we obtain
	\begin{equation*}
		|\Delta w(x)-H(x)|\leq C|x|^{2s-n},\ \ \ |x|>2R_1.
	\end{equation*}
	Therefore, \eqref{eqHl} yields
	\begin{equation*}
		|\Delta w(x)|\leq C|x|^{-2\ep},\ \ \ |x|>2R_1.
	\end{equation*} 
	Since $2\ep<1$, Remark \ref{remdec} implies that 
	\begin{equation*}
		|D^kw(x)|\leq C|x|^{2-k-2\ep},\ \ \ |x|>2R_1,\ \ \ k=0,1,2.
	\end{equation*}
	Finally we can replace \eqref{eqv1} in Lemma \ref{lem21} by the condition that $|w(x)|\leq C|x|^{2-2\ep}$ for $|x|>2R_1$ 
	to complete the proof of this lemma. 
	\qed

	Now we can prove Theorem \ref{thm1}.
	
	\noindent{\it Proof of Theorem \ref{thm1}.} The proof was divided into three steps:
	
	{\bf Step 1.} We show that there exists a vector $b\in\rn$ such that 
	\begin{equation}\label{eqb1}
		Dw(x)\rightarrow b\ \ \ \text{as}\ |x|\rightarrow+\infty.
	\end{equation}
	
	By choosing $\ep$ small we can find a positive integer $m_0$ such that 
	\begin{equation*}
		2^{m_0}\ep<1<2^{m_0+1}\ep.
	\end{equation*}
	Let $\ep_1=2^{m_0}\ep$. It is obvious that $1<2\ep_1<2$. One can apply Lemma \ref{lem22} $m_0$ times to obtain
	\begin{equation}\label{eqpw1}
		\left\{
		\begin{aligned}
			&|D^kw(x)|\leq C|x|^{2-k-\ep_1},\ \ \ k=0,1,2,\ \ \ |x|>2R_1,\\
			&\frac{|D^2w(x_1)-D^2w(x_2)|}{|x_1-x_2|^{\alpha}}\leq C|x_1|^{-\ep_1-\alpha},\ \ \ |x_1|>2R_1,
			\ \ |x_1-x_2|\leq\frac{|x_1|}{2}.
		\end{aligned}
		\right.
	\end{equation}
	Let $F$ and $H$ be as in Lemma \ref{lem22}, then, similar to \eqref{eqF}-\eqref{eqHl}, we have
	\begin{gather*}
		|F(x)|\leq C|x|^{-2s-2\ep_1},\ \ \ |x|>2R_1,\\
		|H(x)|\leq C|x|^{-2\ep_1},\ \ \ |x|>2R_1.
	\end{gather*}
	Note that $1<2\ep_1<2$ and $2s-n<-2$, we can apply the comparing argument as in Lemma \ref{lem22} to obtain 
	\begin{equation*}
		|\Delta w(x)|\leq|H(x)|+C|x|^{2s-n}\leq C|x|^{-2\ep_1}
	\end{equation*}
	and thus $|D^2w(x)|\leq C|x|^{-2\ep_1}$. Since $2\ep_1>1$, \eqref{eqb1} follows from Lemma \ref{lemclaim}.
	
	{\bf Step 2.} Define $w_1(x)=w(x)-b\cdot x$. We will find a $c\in\rr$ such that
	\begin{equation*}
		w_1(x)\rightarrow c\ \ \ \text{as}\ |x|\rightarrow+\infty.
	\end{equation*} 
	
	Note that the equation for $w_1$ can be written as 
	\begin{equation}\label{eqpweq1}
		\tilde{a}_{ij}\partial_{ij}w_1=f-1.
	\end{equation}
	Since $\partial_{ij}w_1\in C^\alpha(\rn)$, we can apply $(-\Delta)^s$ to both side of \eqref{eqpweq1} to get  
	\begin{equation}\label{eqpweqs}
		(-\Delta)^s\left(\Delta w_1\right)=F_1
	\end{equation}
	with 
	\begin{equation*}
		\begin{aligned}
			F_1(x):=&(-\Delta)^sf(x)-\partial_{ij}w_1(x)(-\Delta)^s\tilde{a}_{ij}(x)
			-(\tilde{a}_{ij}(x)-\delta_{ij})(-\Delta)^s(\partial_{ij}w_1)(x)\\
			&+\int_{\rn}\frac{\left(\tilde{a}_{ij}(x)-\tilde{a}_{ij}(y)\right)
				\left(\partial_{ij}w_1(x)-\partial_{ij}w_1(y)\right)}{|x-y|^{n+2s}}dy.
		\end{aligned}
	\end{equation*}
	It follows from \eqref{eqpw1} that 
	\begin{equation*}
		|F_1(x)|\leq C|x|^{-2s-2\ep_1},\ \ \ |x|>2R_1.
	\end{equation*}
	Now we construct $H_1$ as that of $H$ such that $(-\Delta)^sH_1=F_1$. Then 
	\begin{equation*}
		|H_1(x)|\leq C|x|^{-2\ep_1},\ \ \ |x|>2R_1.
	\end{equation*}
	Since $\Delta w_1-H_1\rightarrow0$ near infinity, we obtain by the comparing argument that 
	\begin{equation}\label{eqpwH1}
		|\Delta w_1(x)-H_1(x)|\leq C|x|^{2s-n},\ \ \ |x|>2R_1.
	\end{equation}
	It follows $|\Delta w_1(x)|\leq C|x|^{-2\ep_1}$ for $|x|>R_1$ and hence $|D^2w_1(x)|\leq C|x|^{-2\ep_1}$.
	Going back to \eqref{eqpweqs}, the new estimates for $F_1$ are 
	\begin{equation*}
		|F_1(x)|\leq C|x|^{-2s-\beta}+C|x|^{-2s-4\ep_1},\ \ \ |x|>2R_1.
	\end{equation*}
	Thus, $H_1$ admits a faster decay
	\begin{equation*}
		|H_1(x)|\leq C\left(|x|^{-\beta}+|x|^{-4\ep_1}\right),\ \ \ |x|>2R_1.
	\end{equation*}
	Since $n\geq3$, combined with \eqref{eqpwH1} we obtain
	\begin{equation}\label{eqg1}
		|\Delta w_1(x)|\leq C\left(|x|^{2s-n}+|x|^{-\beta}+|x|^{-4\ep_1}\right)\leq C|x|^{-2-\nu}
	\end{equation}
	for some $\nu>0$. Note that $Dw_1(x)\rightarrow0$ near infinity, we can use Lemma \ref{lemclaim} twice to
	find the expected $c$. 
	
	{\bf Step 3.} Let $w_2(x):=w_1(x)-c=w(x)-b\cdot x-c$. We show 
	\begin{equation*}\label{eqstep3}
		|w_2(x)|\leq 
		\left\{
		\begin{aligned}
			C|x|^{2-\min\{\beta,n\}},\ \ \ \beta\neq n,\\
			C|x|^{2-n}(\ln|x|),\ \ \ \beta=n,
		\end{aligned}
		\right.\ \ \ \ \ \text{for}\ |x|>2R_1.
	\end{equation*}
	
	One can easily check that $|w_2(x)|\leq C$ for $|x|>2R_1$ and $w_2(x)\rightarrow0$ near infinity. Combining with Lemma \ref{lem21}, 
	Lemma \ref{lemclaim} and \eqref{eqg1} we obtain 
	\begin{equation*}
		|D^kw_2(x)|\leq C|x|^{-k-\nu},\ \ \ k=0,1,2,\ \ \ |x|>2R_1.
	\end{equation*}
	The equation for $w_2$ can be written as 
	\begin{equation*}
		\tilde{a}_{ij}\partial_{ij}w_2=f-1,
	\end{equation*}
	where $\tilde{a}_{ij}$ satisfies 
	\begin{align*}
		|\tilde{a}_{ij}(x)-\delta_{ij}|\leq C|x|^{-2-\nu},\ \ \ |x|>2R_1.
	\end{align*}
	We rewrite this equation as 
	\begin{equation*}
		\Delta w_2=F_2:=f-1-(\tilde{a}_{ij}-\delta_{ij})\partial_{ij}w_2.
	\end{equation*}
	It is easy to check that
	\begin{equation*}
		|F_2(x)|\leq C|x|^{-\beta}+C|x|^{-4-2\nu},\ \ \ |x|>2R_1.
	\end{equation*}
	Define 
	\begin{equation*}
		H_2(x):=-c_{n}\int_{\rn\setminus B_{R_1}}\frac{F_2(y)}{|x-y|^{n-2}}dy.
	\end{equation*}
	Then $\Delta H_2(x)=F_2(x)$ in $|x|>2R_1$ and 
	\begin{equation*}\label{eqg2}
		|H_2(x)|\leq
		\left\{
		\begin{aligned}
			C|x|^{2-\min\{\beta,n\}}+|x|^{-2-\nu},\ \ \ \beta\neq n,\\
			C|x|^{2-n}(\ln|x|)+|x|^{-2-\nu},\ \ \ \beta=n.
		\end{aligned}
		\right.
	\end{equation*}
	This kind of decay can be obtained from the $s=1$ version of Corollary \ref{cor2}, i.e., 
	the decay for Newton potential. We remind the reader to see \cite[Lemma 3.2]{lb}. 
	A similar comparing argument as before shows that
	\begin{equation*}
		|w_2(x)-H_2(x)|\leq C|x|^{2-n},\ \ \ |x|>2R_1,
	\end{equation*}
	which gives that 
	\begin{equation*}
		|w_2(x)|\leq 
		\left\{
		\begin{aligned}
			C\left(|x|^{2-\min\{\beta,n\}}+|x|^{-2-\nu}\right),\ \ \ \beta\neq n,\\
			C\left(|x|^{2-n}(\ln|x|)+|x|^{-2-\nu}\right),\ \ \ \beta=n,
		\end{aligned}
		\right.\ \ \ \ \ \text{for}\ |x|>2R_1.
	\end{equation*}
	If $|x|^{-2-\nu}>|x|^{2-n}+|x|^{2-\beta}$ (or $|x|^{-2-\nu}>|x|^{2-n}(\ln|x|)$ ), we can repeat {\bf Step 3} finite 
	times to remove the $|x|^{-2-\nu}$ from the above estimate. 
	
	Eventually by Lemma \ref{lem21} we can establish Theorem \ref{thm1}. 
	\qed
	
	%


	\subsection*{Acknowledgments}
	The authors thank the Referees for their very valuable comments.
	
	This work was supported by 
	the National Natural Science Foundation of China (12371200).

	\noindent\textbf{Data availability} 
	
	\noindent Data sharing not applicable to this article as no datasets were generated or analysed for this
	manuscript.
	
	\noindent\textbf{Declarations}
	
	\noindent{\bf Conflict of interest} The authors have no relevant financial or non-financial interests to disclose.


	%

\end{document}